\documentclass[12pt, oneside, a4paper]{article}
\usepackage[english]{babel}  
\usepackage{float}
\usepackage[utf8]{inputenc}  
\usepackage[T2A]{fontenc}      
\usepackage{wrapfig}
\usepackage{amsmath, amssymb, amsthm, amsfonts, amsthm, dsfont,mathrsfs}
\usepackage{graphicx} 
\usepackage[]{hyperref}
\usepackage{epigraph}
\usepackage{accents}
\usepackage{abstract}
\DeclareMathOperator{\dist}{dist}

\DeclareMathOperator{\grad}{grad}

\newcommand{\eps}{\varepsilon}
\newcommand{\R}{{\mathbb R}}

\newtheorem{theorem}{Theorem} 
\newtheorem{lemma}[theorem]{Lemma}
\newtheorem{define}[theorem]{Definition}
\newtheorem{predl}[theorem]{Proposition}

\usepackage{fancyhdr}
\makeatletter
\def\blindfootnote{\gdef\@thefnmark{}\@footnotetext}
\makeatother
\usepackage{mathtools}

\begin{document}
	\thispagestyle{empty}
	
	\title{On spectra of hyperbolic surfaces without thin handles }
	\author{M.B. Dubashinskiy\footnote{Chebyshev Laboratory, St.~Petersburg State University, 14th Line 29b, Vasilyevsky Island, Saint~Petersburg 199178, Russia.}} 
	\date{\today}
	\maketitle

\renewcommand{\abstractname}{}

\begin{abstract}
	\noindent {\bf Abstract.} We obtain a sharp lower estimate on eigenvalues of Laplace--Beltrami operator on a hyperbolic surface with injectivity radius bounded from the below.
	
\end{abstract}

\thispagestyle{empty}
\blindfootnote{\hspace{0.2mm}e-mail: \texttt{mikhail.dubashinskiy@gmail.com}}
\blindfootnote{Research is supported by the Russian Science Foundation grant 14-21-00035.}
\blindfootnote{{Keywords:} \emph{hyperbolic surface, Laplace--Beltrami eigenvalues, Cheeger-Yau isoperimetric inequality}.}
\blindfootnote{\hspace{0.0mm}{MSC 2010}: primary 35P15, secondary 58J50.}
	
\section{Introduction}

\noindent Let $\Omega$ be a hyperbolic surface, that is, a Riemannian manifold of real dimension $2$ with constant Gaussian curvature $-1$; we assume that $\Omega$ is compact and has no border. Denote by $g$ the genus of $\Omega$. Let $\Delta$ be Laplace--Beltrami operator on $\Omega$; it has purely discrete spectrum since $\Omega$ is compact. Denote by $\lambda_j=\lambda_j(\Omega)$ the $j$th eigenvalue of $-\Delta$ ($j=0,1,2,\dots$). Our main result is the following
 
\begin{theorem}
	\label{th:main}
	Let $r>0$.	There exists a constant $c(r)>0$ 	such that if injectivity radius of $\Omega$ is greater than $r$ then 	$\lambda_{\lceil{\eps g}\rceil} \ge c(r)\cdot \eps^2$ for any $\eps \le 2$.
\end{theorem}

\noindent In what follows, we denote by $c(r)$ any positive constant depending only on $r$ (but not on $\eps$, $g$ and $\Omega$).

Proposition \ref{prop:sharp} below shows that our estimate is sharp in the order.

A theorem by Otal and Rosas (\cite{OtalRosas}) says that $\lambda_{2g-2} > 1/4$ for any $\Omega$ of genus $g$. To the other hand, for a given $\delta >0$, $N\in \mathbb N$ and $g=2,3,\dots$ there exists a hyperbolic surface $\Omega$ of genus $g$ with $\lambda_{2g-3} < \delta$ and $\lambda_{2g-2+N} < 1/4+\delta$. Validity of these inequalities is related to the existence of thin handles on $\Omega$ (see \cite{Buser}). In other words, eigenvalues are small when injectivity radius of $\Omega$ degenerates.  Theorem \ref{th:main}  gives the lower estimate on  eigenvalues under the assumption on this radius. 

\medskip

\noindent {\bf Acknowledgments.} The problem was stated by M. Mirzakhani. Author is also grateful to P.G. Zograf for introducing the topic.

\section{Proof of Theorem \ref{th:main}}

Our proof of Theorem \ref{th:main} is a slight refinement of Buser's argument leading to the estimate $\lambda_{2g-2} \ge 10^{-12}$ (see \cite{Buser}) together with simple Lemma \ref{lemma:graph} on graphs. 

We are going to apply Dirichlet--Neumann bracketing technique. Recall that if $X\subset\Omega$ is a set with positive area and piecewise smooth boundary, then its \emph{Cheeger constant} is defined as
$$
h(X) := \inf\frac{l(A)}{\min\{|B|, |B'|\}},
$$
where $A$ ranges over the family of all finite unions of piecewise smooth curves on $X$ cutting $X$ into two disjoint subsets $B$ and $B'$. Here, $l(A)$ is length of $A$ and $|\cdot|$ is Riemannian volume on $\Omega$. A very standard combination of geometric implementation of minimax principle together with Cheeger--Yau isoperimetric inequality (\cite{Buser}, see also \cite{Cheeger}, \cite{Yau}) leads to the following conclusion: 

\begin{theorem}
	\label{th:Cheeger_minimax}
	Suppose that $k\in\mathbb N$ and that $\Omega$ is subdivided into union of sets $X_1, \dots, X_k$ with piecewise-smooth boundaries and disjoint interiors. Then 	
	$$
	\lambda_k(\Omega) \ge \min\limits_{j=1,\dots, k} \frac{h^2(X_j)}4.
	$$
\end{theorem}

\noindent An appropriate  subdivision of $\Omega$ will be obtained via trianguation of controlled size. For this, recall a result by Buser (Theorem 4.5.2 in \cite{Buser}, see also \cite{Buser78}).

\begin{define}
	\label{def:trigon}
	A closed domain $D\subset \Omega$ is called a \emph{trigon} if $D$ is of one of the following two types: 
	\begin{enumerate}
		\item $D$ a simply connected embedded geodesic triangle \emph{(an ordinary triangle)};
		\item $D$ is embedded doubly connected domain bounded by a geodesic cycle and two geodesic arcs \emph{(a collar-type trigon)}.
	\end{enumerate} 
	Geodesic boundary components of such $D$ are called \emph{sides} of $D$.
\end{define}

\begin{theorem}[Buser]
	\label{th:triangulation}
	Surface $\Omega$ can be triangulated into trigons having side lengths $\le 	\log 4$ and areas between $0.19$ and $1.36$. 
\end{theorem}

\noindent Fix such a triangulation; denote by $\mathcal T_c$ and $\mathcal T_t$ the sets of its collar-type trigons and ordinary triangles respectively. Also, denote by $\mathcal S_c$ and $\mathcal S_a$ the sets of sides of our triangulation which are cycles and geodesic arcs respectively. Let $\mathcal N$ be the set of vertices of triangulation. The proof of Theorem \ref{th:triangulation} from~\cite{Buser} furnishes symmetries of trigons from $\mathcal T_c$: namely, $\mathcal S_a$-sides of such a trigon have equal lengths. From this we derive that lengths of arcs from $\mathcal S_a$ are bounded from the below by an absolute constant; also, angles of triangulation are also bounded from below by an absolute constant. (For $\mathcal T_t$-trigons these statements are obvious due to upper area estimate whereas for segments and angles in boundaries of collars the computation is done in \cite{Buser}.)

\begin{lemma}
	\label{lemma:shortcut}
	If $a_1,a_2$ are two sides of triangulation with no common vertex, then $\dist_\Omega(a_1, a_2)$ is bounded from the below by an absolute constant $d_0>0$. 
\end{lemma}

\noindent {\bf Proof.} First, notice that distance from any $c\in \mathcal S_c$ to any other side is bounded from below by a universal constant --- otherwise area of some trigon from $\mathcal T_c$ degenerates.

Next, we claim that \emph{distances between vertices of triangulation are bounded from below by a universal constant}. Indeed, let $U$ be a metric ball on $\Omega$ centered in some $v\in\mathcal N$. If radius of $U$ is small enough then for each $\tau \in \mathcal T_c \cup \mathcal T_t$ intersection $U\cap \tau$ can intersect no sides of triangulation except for those who emanate from $v$; it is easily checked for both types of trigons, and this leads to our claim.

Now, suppose that $\gamma$ is a geodesic arc joining $a_1$ and $a_2$ and of small length; it cannot intersect some side from $\mathcal S_c$ since such sides are far away enough from all the other sides. Suppose that $\gamma$ passes through some trigon $\tau\in\mathcal T_c\cup \mathcal T_t$. Then it occurs close enough to some vertex $v\in\mathcal N$ (because angles of trigons are bounded from below). Since vertices are separated, the whole curve $\gamma$ is situated close enough to some vertex $v\in\mathcal N$, but in this case $\gamma$ can join only sides emanating from $v$. Proof is finished.
$\blacksquare$

\medskip

\noindent Now we estimate Cheeger constants:

\begin{lemma}
	\label{lemma:Cheeger}
	Let $X\subset\Omega$ be a union of $N$ distinct trigons from our triangulation \emph{(}$N=1,2,\dots$\emph{)}. Suppose that $X$ is "connected" in the sense that two trigons are adjacent if they have a common side, not just a vertex. \emph{(}More formally, we may say that the interior of $X$ is connected.\emph{)}
	
	Then, under hypothesis on injectivity radius of $\Omega$, we have 
   \begin{equation}
    \label{eq:isoper}
	h(X) \ge \frac{c(r)}N.
	\end{equation}
\end{lemma}

\noindent {\bf Proof.} 
Let $A, B, B'$ be sets from definition of Cheeger constant for $X$; we have $A\neq \varnothing$ since $X$ is connected. By Yau lemma (\cite{Buser}, Lemma 8.3.6, see also \cite{Yau}) we may assume that $B, B'$ are connected. If $l(A) \ge r$ then note that $\min\{|B|, |B'|\} \le 1.36\cdot 1/2\cdot N$, and this leads to~(\ref{eq:isoper}). Next, suppose that $A$ contains a cycle $\gamma$. Then  $\gamma$ is homotopic to identity in $\Omega$ (since  $l(A) < r$ and by injectivity radius condition). Cycle $\gamma$ should enclose in $\Omega$ a component of area $\le l(\gamma) / h(\mathds H)=l(\gamma)$ (it is known that Cheeger constant of the whole Lobachevskiy plane $\mathbb H$ is $1$)  and this also gives (\ref{eq:isoper}). So, suppose that $A$ does not contain a cycle.

We could assume from the beginning that $r < d_0$ where $d_0$ is the constant from Lemma~\ref{lemma:shortcut}. Set $A$ is a union of curves; take any component $\gamma$ of $A$. Then $\gamma$ necessarily has ends (since $A$ does not contain a cycle) and these ends lie on $\partial X$. Take two of such ends, $p_1, p_2$, and curve $\gamma_1\subset \gamma$ joining them. By Lemma~\ref{lemma:shortcut}, $p_1$ and $p_2$ are situated either on the same side of triangulation or on two distinct sides emanating from their common vertex; this side or these sides lie on $\partial X$. But if $Y$ is an angle on $\mathbb H$ or half-plane of $\mathbb H$ then $h(Y)=1$ (see, e.g., proof of Theorem 8.1.2 in \cite{Buser}). This and also injectivity radius condition, say, in $p_1$ lead to~(\ref{eq:isoper}). $\blacksquare$

\medskip

\noindent Now, to obtain a subdivision of $\Omega$ via our triangulation, we give a simple graph lemma:

\begin{lemma}
	\label{lemma:graph}
	Let $G$ be a finite connected non-oriented graph with degrees of vertices  $\le 3$. Let $k \in\mathbb N$. 
	The set of vertices of $G$ can be subdivided as $V_1 \sqcup V_2 \sqcup \dots \sqcup V_\alpha \sqcup V'$ \emph{(}for some $\alpha=0,1,2,\dots$\emph{)} such that:
	\begin{enumerate}
		\item graphs induced by $G$ on each $V_1, V_2, \dots, V_\alpha, V'$ are connected;
		\item $2^k \le |V_1|, |V_2|, \dots, |V_\alpha| \le 2^{k+1}-1$ and $0 \le |V'| \le 2^k$.
	\end{enumerate}
\end{lemma}
	
\medskip

\noindent {\bf Proof.} We argue by induction by the number of vertices in $G$; for the empty graph the statement is obvious. We may assume that $G$ is a tree. Pick a leaf of $G$ and call it root. Arrange the graph by levels by distance from the root. Vertex $v$ from some level is adjacent to $\le 2$ vertices from the next level, we call them \emph{children} of $v$.

Let us construct a sequence of vertices $v_0, v_1, \dots, v_\beta$ of $G$ ($\beta$ will be some non-negative integer). Take the root as $v_0$. Suppose that $v_j$ is constructed and that $v_l$ and $v_r$ are its children. W.l.o.g., the total number of descendants of $v_l$ is greater or equal than that of $v_r$. Then put $v_{j+1}:=v_l$. If $v_j$ has only one child then take it as $v_{j+1}$; and if $v_j$ has no children then stop our process and put $\beta:=j$, this should occur necessarily. Thus we  construct a sequence of vertices.

Now pass this sequence in the reverse order (starting from $v_\beta$ and up to $v_0$) and watch for the total number of descendants of vertices. If $v_{j+1}$ has $x$ descendants (together with itself) then $v_{j}$ has $\le 2x+1$ descendants together with itself. Then we have two cases:

\begin{enumerate}
	\item There exists some $v_j$ having $\ge 2^k$ and $\le 2^{k+1}-1$ descendants together with itself. Then, for $V_1$ we take the set consisting of $v_j$ and of all of its descendants. Cut them from $G$ and apply induction hypothesis for $G$ without $V_1$.
	
	\item $|G|<2^k$. Then take $V'$ as the whole set of vertices of $G$. $\blacksquare$
\end{enumerate}

\noindent {\bf Proof of Theorem \ref{th:main}.} First, assume that $\eps g\le 1$. Then we have to prove that $\lambda_1 > c(r)/g^2$, but, by Theorem \ref{th:Cheeger_minimax}, it is enough to prove that $h(\Omega) > c(r)/g$. Taking $A$ from the definition of $h(\Omega)$, we see that $A$ must contain a cycle; in this case we argue as in the corresponding case in the proof of Lemma \ref{lemma:Cheeger} and easily obtain the desired (recall that $|\Omega|=2\pi(2g-2)$).

Now, suppose that $\eps g >1$. Pick $k\in\mathbb N$ with $2^k\ge \dfrac{8\pi}{0.19\cdot \eps} \ge 2^{k-1}$, this can be done because $\eps\le 2$. Let $G$ be the graph of triangulation obtained in Theorem~\ref{th:triangulation}: namely, set of vertices of $G$ is $\mathcal T_t\cup \mathcal T_c$ and two such trigons are adjacent if they have a common side. Apply Lemma~\ref{lemma:graph} to $G$, take subdivision of the set of vertices of $G$ obtained by this lemma and consider corresponding subdivision of $\Omega$ as $X_1 \cup X_2 \cup \dots \cup X_\alpha \cup X'$ for some $\alpha =0,1,2,\dots$. Since trigons have area $\ge 0.19$, we have $|X_j| \ge 2^k\cdot 0.19$ for all $j=1,2,\dots,\alpha$. Then 
$$\alpha \le \dfrac{|\Omega|}{2^k\cdot 0.19} < \dfrac{4\pi g}{2^k\cdot 0.19} \le \frac{\eps g}{2} \le \lceil \eps g\rceil-1.$$ So, $\alpha+1\le  \lceil \eps g\rceil$. 
Now, by Lemma \ref{lemma:Cheeger}, we have $h(X_j), h(X') \ge c(r)/2^k$ for all $j$. This and Theorem \ref{th:Cheeger_minimax} lead to the desired. $\blacksquare$

\medskip

\medskip

\noindent Finally, let us demonstrate the sharpness of our estimate (we may think  that $\eps$ is $1/k$).

\begin{predl}
	\label{prop:sharp}
	For any $k, l \in \mathbb N$ there exists a hyperbolic surface $\Omega$ of genus $kl+1$ with injectivity radius bounded from below by a universal constant and with $\lambda_{l-1}(\Omega) \le C/k^2$, here $C<+\infty$ is a universal constant.
\end{predl}
	
\noindent {\bf Proof.} Let $P$ be fixed hyperbolic pants bounded by geodesic cycles of length, say, $1$ (existence and uniqueness of such pants is a well-known fact). Let $P_1, \dots, P_{2k}$ be copies of these pants. For $j=1,2, \dots, 2k$, denote by $\gamma_1(P_j), \gamma_2(P_j), \gamma_3(P_j)$ the boundary components of~$P_j$. For $j=1,2,\dots, k$, let us glue $\gamma_2(P_{2j-1})$ to $\gamma_2(P_{2j})$ and  $\gamma_3(P_{2j-1})$ to $\gamma_3(P_{2j})$. Also, for $j=1,2,\dots, k-1$ paste $\gamma_1(P_{2j})$ to $\gamma_1(P_{2j+1})$. Denote by $Q$ the surface obtained in such a way; it is a hyperbolic surface with two geodesic boundary components of length~$1$. There exists a Sobolev function $f\colon Q\to \R$ with the following properties: first, $f=0$ on $\partial Q$; second, $f$ takes values in $[j-1,j]$ on ${P_j}$ and on $P_{2k-j+1}$ for $j=1,2, \dots, k$; third, $|\grad f|$ does not exceed some absolute constant, $\grad$ being metric gradient. (To construct such a function, just let it be equal to appropriate constants on boundary components of pants and interpolate it into the interiors of pants anyway.) Now, take $l$ copies of $Q$ and paste them in a cyclic way to obtain a hyperbolic surface $\Omega$ with no boundary. Then genus of~$\Omega$ is $kl+1$. Moreover, one can find Sobolev functions $f_1, f_2, \dots, f_l\colon \Omega\to \R$ with disjoint supports and such that $\int_\Omega f_j^2 \ge c_1\cdot k^3$, $\int_\Omega |\grad f_j|^2 \le c_2\cdot k$ (constants $c_1, c_2$ are absolute). By the geometric version of minimax principle (that is, by upper estimate from Dirichlet--Neumann bracketing), this leads to the desired eigenvalue estimate. Injectivity radius of $\Omega$ is bounded from the below since it is true for any pants. $\blacksquare$

\end{document}